\documentclass[11pt]{amsart}

\usepackage{amsmath,amssymb}

\newtheorem{Lemma}{Lemma}[section]

\newtheorem{theorem}{Theorem}[section]
\newtheorem{prop}[theorem]{Proposition}

\newtheorem{remark}[theorem]{Remark}
\newtheorem{example}[theorem]{Example}

\def\beq{\begin{equation}}
\def\eeq{\end{equation}}
\def\beqn{\begin{eqnarray}}
\def\eeqn{\end{eqnarray}}
\def\bear{\begin{eqnarray*}}
\def\eear{\end{eqnarray*}}

\def\bl{\begin{trivlist}}
\def\el{\end{trivlist}}

\def\bm{\begin{pmatrix}}
\def\em{\end{pmatrix}}

\def\noi{\noindent}

\def\al{\alpha}
\def\be{\beta}

\def\phi{\varphi}

\def\R{{\mathbb R}}

\def\rg{\rangle}

\title{Interacting two-state Markov chains on undirected networks }
\author{Maung Min-Oo}
\address{Department of Mathematics \& Statistics,  McMaster University}
\email{minoo@mcmaster.ca}

\date{\today}

\begin{document}

\maketitle

\begin{abstract}
\noi
{\it It is shown that irreducible two-state continuous-time Markov chains interacting on a network in a bilinear fashion have a unique stable steady state. The proof is elementary and uses the relative entropy function.}
\end{abstract}

\bigskip

\section{{\bf Description of the main result }}

\medskip

This is an elementary paper about two-state Markov chain attached to each node (vertex) of a finite undirected network (simple weighted undirected graph). In this paper, we will deal with Markov chains in continuous time (sometimes called Markov jump processes). The interaction between two chains that are linked by an edge of the network is a simple bilinear function of the two opposite states (a coupling constant times the product of the probabilities of the two opposite states). The purpose of this paper is to give an elementary proof of the existence and uniqueness of the steady state (or equilibrium). In short, we will give a simple proof of the following fact:

\smallskip
\noi
{\it There exists a unique steady state of irreducible two-state Markov chains which are linked on an undirected network through an interaction term that depends bilinearly on neighbouring opposite states.} \\
\noi
In more precise technical terms:

\medskip
\noi
{\bf Theorem 3.1} 
{\it Let $\alpha, \beta \in \R^N_+$ ,  $\gamma_{0 1} , \gamma_{ 1 0} \geq 0$ and $W$ be a symmetric $N \times N$ matrix with non-negative entries and 
zeros on the diagonal. Then the system of differential equations:
$$\frac{dp^i}{dt}  =  - \, \alpha^i p^i + \beta^i q^i - \gamma_{01} \, p^i \sum_j W_j^i q^j + \gamma_{10} \, q^i \sum_j W_j^i p^j \qquad i = 1, \ldots , N
$$
where $q^i = 1 - p^i$, leaves the $N$-dimensional unit cube $[0 , 1]^N$ invariant and possesses a unique globally stable steady state (equilibrium point) in the interior of $I^N$.}

\bigskip
\noi
This is proved in section 3, after setting up the notation in the next section. In the final section we make some simple remarks about the steady state distribution and discuss some special cases. In the next paper we plan to deal with the case of directed networks (which is the more interesting case).

\bigskip
\bigskip

\section{ \bf{The structure of the equations}}

\medskip

\subsection{Basic notation and terminology}

\medskip

\subsubsection{Markov chains}
There is an extensive theory of Markov chains. Here are two introductory textbooks: \cite{la}\cite{no}.  We will briefly describe what we need to know. 
In continuous time, a time-homogeneous two-state Markov chain (or a Markov jump process) is completely determined by a $2 \times $2 matrix (the infinitesimal transition probability matrix between the two states) :
$$ Q = 
\bm
-\alpha & \alpha \\
\beta & -\beta
\em
$$
where $\alpha \geq 0 $ and $\beta \geq 0$. The time evolution of the probabilities $p$ and $q = 1 - p$ at the two states $|0 \rangle$ and $|1 \rangle$ is then determined by solving the linear differential equation (with constant coefficients):
$$
\frac{d}{dt}(p,q)  = (p,q) Q 
$$
whose solution is simply $(p(t), q(t)) = (p(0), q(0)) e^{t Q}$, and as
$t \rightarrow \infty$, this converges to the steady state: $ ({\bar p}, {\bar q}) = \frac{1}{\alpha + \beta} ( \beta , \alpha ) $. The Markov chain is irreducible and aperiodic provided $\al$ and $\be$ are strictly positive.

\bigskip
\noi
An important function on the one dimensional simplex ($p+q =1$) is the relative entropy function. With respect to the steady state distribution, it is defined as:
$$ E_{\bar{p}}(p) = E_{(\bar{p}, \bar{q})}(p,q) =  - \bar{p} \log \frac{p}{\bar{p}} - \bar{q} \log \frac{q}{\bar{q}} $$
This is also known as the Kullback-Leibler ``distance" (although it is not symmetric and does not satisfy the triangle inequality). 

\medskip
\noi
As a warm-up exercise, let us compute the evolution of this entropy function along the flow: 

\medskip
\noi
 $$
\frac{dE}{dt}  =  (\frac{\bar{q}}{q} - \frac{\bar{p}}{p}) \dot{p} =  - \frac{(- \alpha p + \beta q)^2}{(\alpha + \beta) p q} = 
- \frac{\alpha + \beta}{p q}  (p - \bar{p})^2$$
which is strictly negative unless $p = \bar{p}$. This proves the uniqueness and global stability of the steady state.

\medskip
\noi
The main purpose of this paper is to do a similar calculation on two-state Markov chains interacting bilinearly on a network.

\bigskip

\subsubsection{ Networks and Graphs }
A finite network (or a graph) is a collection of vertices (or nodes), denoted by $\mathcal{V} = \{ v_1, \ldots v_N \}$,  together with a collection of edges, denoted by $\mathcal{E}$, where each edge joins two vertices. For an undirected network we think of each edge $e \in \mathcal{E} $ as an unordered pair of vertices. An edge connecting a vertex to itself is called a loop. In this paper we will consider undirected finite graphs without loops. We will also assume that there is at most one edge between two different vertices, but we will consider the case where each edge is assigned a weight $w(e) = w_{i j} = W^i_j = W^j_i$, a positive real number. We will set $w_{i j} = 0$ if there is no edge between $v_i$ and $v_j$. In particular $w_{i i} = 0$. By default, if there is no specific weight attached, each edge will have weight $1$ (we then sometimes say it is a graph!). The whole information about an undirected network is therefore completely encoded by a real symmetric matrix $ W^i_j$ with non-negative entries and zeros on the diagonal. The sum of $i^{th}$-row of $W$ (which is the same as the sum of the $i^{th}$-column), denoted by $d^i$, is the (weighted) number of edges that contain the vertex $v_i$ and is called the degree of that vertex. We will denote the diagonal matrix of these degrees $d^1, \ldots, d^N$ by $D$. 
 The (combinatorial) Laplacian of the network is now defined as: $L = D - W$. $L$ is a symmetric matrix with non-positive non-diagonal entries with all row sums (and column sums) equal to zero. 
 $L$ determines $W$, and so encodes all the information about the graph. The (non-negative) quadratic form associated to the Laplacian $L$ is then:
$$ < Lx , x > \, = \sum_{i,j} w_{i j} |x^i - x^j|^2
$$
where $x^i = x(v_i) \; i = 1, \ldots, N$ is a function defined on the vertices, thought of as a column vector. Since $L$ vanishes on constant functions, $0$ is always an eigenvalue and it is a simple eigenvalue iff the network is connected. We note that $Lx(v_{max}) \geq 0$ and $ Lx(v_{min}) \leq 0$ if $v_{max}$ and $v_{min}$ are respectively, local maximum and minimum points of $x$. The inequalities are strict for strict maxima and minima.

\medskip
\noi
These are very basic elementary facts about networks and graphs. There is an extensive theory and here are two introductory books: \cite{ch}\cite{ne}.

\bigskip

\subsection {Interacting Markov processes on a network}

Now suppose that at each node (vertex) $v_i$ of a network with weight matrix $W_j^i$, we have a continuos time Markov chain with two states 
$|0 \rangle \, , \, |1 \rangle$ and infinitesimal transition matrix: 
$Q^i = \bm
-\alpha^i & \alpha^i \\
\beta^i  & -\beta^i
\em$. \\

\medskip
\noi
We will denote the probabilities at each node by $(p^i, q^i= 1 - p^i)$. 
If there would be no interaction between the Markov processes at different nodes, the infinitesimal transition matrix for the whole system, consisting of $2^N$ states, 
will be the tensor product acting as a derivation $ {\bf Q} = \sum_i id \otimes \cdots \otimes Q^i \otimes \cdots \otimes id $\,, one node at a time. 
This describes a random walk on the hypercube $\{ 0, 1 \} ^N$.  We will add to this tensor product an interaction term depending on the network as follows. 
At each node $v_i$, we change the probabilities of ${\bf Q } $ by adding the following terms which describe a very simple bilinear interaction between the two states at neighbouring sites (only opposite states will interact). 
\bear
 p^i &\mapsto&   - \, \gamma_{01} \, p^i \sum_j W_j^i q^j  + \gamma_{10} \, q^i \sum_j W_j^i p^j \\
q^i &\mapsto&    + \, \gamma_{01} \, p^i \sum_j W_j^i q^j - \gamma_{10} \, q^i \sum_j W_j^i p^j 
\eear
where $\gamma_{01} \geq 0$ and  $\gamma_{10} \geq 0$ are coupling constants (not necessarily equal). This means that we are changing the transitional probabilities 
of the independent tensor product process by a bilinear interaction term that depends on the network and on the coupling constants between the opposite states.

\medskip
\noi
The new system is strictly speaking not a Markov chain on the hypercube with $2^N$ states, but it can be thought of as a ``non-linear Markov process" on the (continuous) space of all probability distributions on the nodes of the network. The state space is therefore $I^N = [0,1]^N$ and we study the following dynamical system of $N$ (independent) differential equations where the non-linearity is of a simple type.

\begin{equation}
\frac{dp^i}{dt}  =  - \frac{dq^i}{dt}  =  - \, \alpha^i p^i + \beta^i q^i - \gamma_{01} \, p^i \sum_j W_j^i q^j + \gamma_{10} \, q^i \sum_j W_j^i p^j 
\end{equation}

\medskip
\noi
In terms of the Laplacian, these equations can also be written as:

\begin{equation}
\frac{dp^i}{dt}  =    - \frac{dq^i}{dt} = - \alpha^i p^i + \beta^i q^i - \hat{\gamma}\, d^i \, p^i \, q^i   
+ \gamma_{01} \, p^i \sum_j L_j^i q^j   - \gamma_{10} \, q^i \sum_j L_j^i p^j  
\end{equation}
where $ d_i = \sum_j w_{i j} $ is the degree of the vertex $v_i$   and $ \hat{\gamma} =  (\gamma_{01}  -  \, \gamma_{10})  $.

\medskip
\noi
Note that $\sum_j L^i _{j} q^j = - \sum_j L^i_j p^j $ and that $ \sum_i \sum_j L^i_j x^j = 0$ for any $x^j$.\\

\noi
Let us denote the vector $Lp = -Lq$ by $ l $ , i.e.,  $l^i = \sum_j L^i_j p^j$.
We then have:
\begin{equation}
\frac{dp^i}{dt}  =   \hat{\gamma}\, d^i \,( p^i)^2  - (\alpha^i  + \beta^i  + \hat{\gamma}\, d^i)  \, p^i    + \beta^i 
 -  \hat{\gamma} \, p^i \, (Lp)^i -  \gamma_{10}  (Lp)^i 
\end{equation} 

\medskip
\noi
The equations decouple on different connected components of the network.
The relation to the bigger system on the ``hypercube" $I^{2^N}$ is given by the embedding 
$\Phi: I^N  \rightarrow  I^{2^N} \; p= (p^i) \mapsto {\bf p} = \Phi(p) $ defined as:
$$
{\bf p}^{\sigma} = p(\sigma) = \prod_i p^i(\sigma^i) 
$$
where for $(\sigma^i) \in \{0 ,1\}^N$ which is a sequence of $0$'s and $1$'s (corners of $I^N$), we define $ p^i(0) = p^i$ and $ p^i(1) = 1-p^i = q^i$. Obviously, 
$\sum_{\sigma} p(\sigma) = \prod_i (p^i + q^i) = 1$.

\medskip
\noi
The process on $I^{2^N}$ is now determined by its ``restriction" to an $N$-dimensional ``quadratic variety".  Note also that the relative entropy (Kullback-Leibler divergence) on the hypercube, restricted to this subvariety is just the sum of the relative entropies at each node. 

\begin{equation}
 E_{{\bf {\bar p}}}({\bf p} ) = \, - \sum_{\sigma} {\bf {\bar p}}_{\sigma} \log \frac{{\bf p}^i_{\sigma}}{ {{\bf \bar p}}_{\sigma}} = \sum_i E_{{\bar p}^i} (p^i)
\end{equation}
\noi
so the embedding $\Phi$ preserves relative entropies.

\bigskip
\bigskip

\section{ {\bf The existence and stability of the steady state} }

\bigskip
\noi
We will prove existence, uniqueness and stability of the steady state of the differential equation   
\bear
\frac{dp^i}{dt}  =  F(p)^i
 \eear 
on the product space $I^{N}$, where $F$ is the vector field:
\beqn
 F(p)^i = 
 \hat{\gamma}\,d^i \,( p^i)^2  - (\alpha^i  + \beta^i  + \hat{\gamma}\,d^i)  \, p^i    + \beta^i 
 -  \hat{\gamma} \, p^i (Lp)^i -  \gamma_{10}  (Lp)^i  
 \eeqn 

\bigskip
\noi
We now check $F$ on the boundary of the cube $I^N$, consisting of $2^N$ faces where one of the $p^i$'s is equal to $0$ or $1$.

\medskip
\noi
When $p^i = 0$, then $(Lp)^i \leq 0$ and hence $F(p)^i = + \beta^i -  \gamma_{10}  \, (Lp)^i > 0 $. 

\smallskip
\noi 
When $p^i = 1$, then $(Lp)^i \geq 0$ and hence $F(p)^i = - \alpha^i -  \gamma_{01}  \, (Lp)^i < 0 $.

\medskip
\noi
So $F$ points inwards and this proves that there exists at least one zero of the vector field (or steady state) in the interior $(0 , 1)^N$.

\bigskip
\noi 
We prove now that there is a unique globally stable steady state by showing that the entropy function with respect to any steady state is strictly decreasing along the flow until it reaches that steady state. To that purpose, we first need to establish a little fact about the Laplacian acting on the unit cube:

\begin{Lemma}
 If $x^i \in (0,1)$ and $y^i \in (0,1)$ for all $i = 1, \dots, N$, then
$$ \sum_i \left( \frac{y^i}{x^i}\sum_j L_j^i x^j + \frac{x^i}{y^i}\sum_j L_j^i y^j \right)  \leq 0$$
and equality holds iff $x^i = y^i$ for all $i$.
\end{Lemma}

\medskip
\noi
{\bf Proof}: Using the definition of the Laplacian:

$$  \sum_{i,j}  (\frac { y^i } {x^i}  L^i_j x ^j   + \frac {x^i} { y^i}  L^i_j  y ^j ) =  
\sum_{i \sim j}  w_{i j}( \frac { y^i } {x^i} -  \frac  {y^j} {x^j} )( x^i -  x^j)  +  w_{i j}(\frac {x^i}{ y ^i} - \frac {x^j}{ y ^j}) ( y^i - y ^j)
$$
where $i \sim j$ means that $i$ and $j$ are connected by an edge with weight $w_{i j} = w_{j i} > 0$. 

\noi
Now
$$
\frac { y^i } {x^i} -  \frac { y^j} {x^j}  = \frac{ y^i  }{x^i x^j} (x^j - x^i) + \frac{1}{x^j} ( y^i - y^j)
$$
$$
 \frac { x^i } { y^i} -  \frac { x^j} { y^j}  = \frac { x^i }{  y^i y^j } ( y^j - y ^i)  + \frac{1} { y^j} ( x^i - x^j)
$$

\noi
and hence
$$ ( \frac { y^i } {x^i} -  \frac  {y^j} {x^j} )( x^i -  x^j)  + (\frac {x^i}{ y ^i} - \frac {x^j}{ y ^j}) ( y^i - y ^j)$$
$$ = - \frac{ y^i  }{x^i x^j} (x^j - x^i) ^2   + (\frac{1}{x^j} +  \frac{1} { y^j}) ( y^i - y^j) ( x^i - x^j) 
-  \frac { x^i }{ y^i y^j } ( y^j - y^i) ^2 $$

\noi
which is a negative definite quadratic form since 
$$ (\frac{1}{x^j} +  \frac{1} { y^j})^2 \geq  \frac{4}{ x^j y^j } 
$$
 \hfill QED \\

\bigskip
\noi 
Any steady state $\bar{p}^i$ satisfies: 
\beqn
 - \alpha^i \bar{p}^i + \beta^i \bar{q}^i - \hat{\gamma}\,d^i \, \bar{p}^i \, \bar{q}^i   -  
 \gamma_{01} \bar{p}^i \, \bar{l}^i -  \gamma_{10} \bar{q}^i \, \bar{l}^i = 0
\eeqn
for each $i$, where $\bar{l}^i = L^i_j \bar{p}^j$.

\bigskip
\noi
Let
\beqn
 E_{\bar{p}}(p) = - \sum_i ( \bar{p}^i \log p^i + \bar{q}^i \log q^i  ) + \sum_i ( \bar{p}^i \log \bar{p}^i + \bar{q}^i \log \bar{q}^i  )
 \eeqn 
which is the sum of all the relative entropies to the steady state at each node.
\bear
\frac{dE }{dt}  & = &   \sum_i(\frac{\bar{q}^i}{q^i} - \frac{\bar{p}^i}{p^i}) \frac{dp^i}{dt}   \\
&=&  \sum_i   (p^i - \bar{p}^i)  \left(   - \frac{\alpha^i} {q^i} + \frac{\beta^i} {p^i} - \hat{\gamma}\,d^i  -  
\gamma_{01} \frac{l^i}{q^i}  -  \gamma_{10} \frac{l^i}{p^i} \right)  \\ 
& =  & \sum_i   (p^i - \bar{p}^i)  \left( \alpha^i ( \frac{1}{\bar{q}^i} -  \frac{1}{q^i}) - \beta^i( \frac{1}{\bar{p}^i} -  \frac{1}{p^i} ) 
  -  \gamma_{01} (\frac{l^i}{q^i} -  \frac{\bar{l}^i}{\bar{q}^i} ) -  \gamma_{10} (\frac{l^i}{p^i}  - \frac{\bar{l}^i}{\bar{p}^i} )\right) 
 \eear
where $l^i = L_j^ip^j ,\,  \bar{l}^i = L^i_j \bar{p}^j$ and we used the steady state equation 3.6.

\noi
Now 
\bear 
 \gamma_{0 1} \sum_i (p^i - \bar{p}^i) (\frac{l^i}{q^i} -  \frac{\bar{l}^i} {\bar{q}^i} ) = \gamma_{0 1} \sum_i (\bar{q}^i - q^i) (\frac{l^i}{q^i} -  \frac{\bar{l}^i} {\bar{q}^i})
 =  - \gamma_{0 1} \sum_i \left( \frac{\bar{q}^i}{q^i} l^i  + \frac{q^i}{\bar{q}^i} \bar{l}^i \right)\\
  \gamma_{1 0} \sum_i (p^i - \bar{p}^i) (\frac{l^i}{p^i} -  \frac{\bar{l}^i} {\bar{p}^i} ) =  \gamma_{1 0} \sum_i (p^i - \bar{p}^i) (\frac{l^i}{p^i} -  \frac{\bar{l}^i} {\bar{p}^i})
  =   \gamma_{1 0} \sum_i \left( \frac{\bar{p}^i}{p^i} l^i  + \frac{p^i}{\bar{p}^i} \bar{l}^i \right)
\eear
since $\sum_i l^i = \sum_i \bar{l}^i = 0$

\medskip
\noi
Using the fact that $l^i = \sum_j L^i_j p^j =  -\sum_j L^i_j q^j \, , \; \, \bar{l}^i = \sum_j L^i_j \bar{p}^j =  -\sum_j L^i_j \bar{q}^j $ we can now apply the basic Lemma above to the terms involving the Laplacian to get:

\bear
\frac{dE }{dt}  & \leq &  \sum_i   (p^i - \bar{p}^i)  \left (  - {\alpha^i} (\frac{1}{q^i} - \frac{1}{\bar{q}^i}) + \beta^i ( \frac{1} {p^i}  - \frac{1}{\bar{p}^i} ) \right) \\
& = &  - \sum_i \, \frac{\alpha^i}{\bar{q}^i q^i} ( \bar{q}^i - q^i)^2  - \sum_i \frac{\beta^i}{\bar{p}^i p^i } (p^i - \bar{p}^i)^2 \\
&\leq& 0    
\eear
with strictly inequality unless  $ p^i = {\bar p}^i $, for all $i$. 

\bigskip
\noi
This proves both uniqueness and global stability of the steady state and hence the following theorem is now established.  

\begin{theorem}
Let $\alpha, \beta \in \R^N_+$ , \, $\gamma_{0 1} , \gamma_{ 1 0} \geq 0$ and $W$ be a symmetric $N \times N$ matrix with all entries non-negative and with 
zeros on the diagonal. Then the system of differential equations:
$$\frac{dp^i}{dt}  =  - \, \alpha^i p^i + \beta^i q^i - \gamma_{01} \, p^i \sum_j W_j^i q^j + \gamma_{10} \, q^i \sum_j W_j^i p^j \qquad i = 1, \ldots , N
$$
where $q^i = 1 - p^i$, leaves the $N$-dimensional unit cube $[0 , 1]^N$ invariant and possesses a unique globally stable steady state (fixed point) in the interior of $I^N$.
\end{theorem}

\bigskip
\bigskip

\section{\bf{Remarks}}

\subsection{The spatial distribution of the steady state}

\medskip
\noi
Let us denote the mean (average) of a function $x$ on the network by $\langle x \rangle = \frac{1}{N}\sum_i x^i$. Let $r = x - \langle x \rangle$. Then the variance of $x$ is given by 
$Var(x) = \langle r^2 \rangle$.
We then have by the basic properties of the Laplacian, the basic inequality:
\beq
\frac{1}{N} \langle x , Lx \rangle = \frac{1}{N}\sum_i x^i \, L_j^i x^j = \frac{1}{N}\sum_i r^i \, L_j^i r^j \geq \lambda_1 Var(x)
\eeq
where $\lambda_1$ is the first positive eigenvalue (which is the same as the second eigenvalue, since we are assuming that the graph is connected) of the Laplacian.

\medskip
\noi
Now since the steady state $\bar{p}$ satisfies:
$$\hat{\gamma}\,d^i\, \bar{p}^i \bar{q}^i   - (\alpha^i + \beta^i ) \, \bar{p}^i  + \beta^i = 
 \hat{\gamma} \bar{p}^i \, L_j^i \bar{p}^j +  \gamma_{10} L_j^i \bar{p}^j$$
 we get by taking averages, the following
 \beq
 - \hat{\gamma} \frac{1}{N}\sum_i d^i \bar{p}^i \bar{q}^i - \frac{1}{N}\sum_i (\alpha^i + \beta^i) \bar{p}^i + \langle \beta \rangle = \hat{\gamma} \frac{1}{N}\sum_i \bar{r}^i \, L_j^i \bar{r}^j
 \eeq
where $ \bar{r} = \bar{p} - \langle \bar{p} \rangle$,  and hence (trivially):
\begin{prop}
The variance of the equilibrium distribution satisfies the estimates:
$$Var(\bar{p}) \leq \frac{1}{\lambda_1} \frac{\langle \beta \rangle}{\hat{\gamma}} \qquad \mbox{if} \; \; \hat{\gamma} > 0 \qquad
\left( Var(\bar{p}) \leq - \frac{1}{\lambda_1} \frac{\langle \alpha \rangle}{\hat{\gamma}} \; \;  \mbox{if} \; \; \hat{\gamma} < 0 \right) $$
\end{prop}

\medskip
\noi
These are not very useful estimate unless $\lambda_1 |\hat{\gamma}|$ is very large.
We will discuss the special case $\hat{\gamma} = 0$ in the next section.\\

\medskip
\noi
Let $R^i$ be the quadratic function: 
\beqn 
R^i (x) = \hat{\gamma}\,d^i\, x^2  - (\alpha^i + \beta^i + \hat{\gamma}\,d^i) \, x  + \beta^i 
\eeqn
defined at each node with a (unique) zero $\rho^i \in (0,1)$. 
Since $L(\bar{p})^i \geq 0 $ at the nodes where $\bar{p}$ attains a local maximum and $L(\bar{p})^i$ at the node where $\bar{p}$ attains a local minimum 
(with strict inequalities for strict (local) maxima and minima nodes we have the following bounds on the absolute maximum and minimum values $\bar{p}_{max}$
and $\bar{p}_{min}$ of the steady state $\bar{p}$. 
\beqn R^{i_{max}}(\bar{p}_{max}) \geq 0  \qquad \mbox{and} \qquad R^{i_{min}}(\bar{p}_{min}) \leq 0 \eeqn

\medskip
\noi
To simplify the discussion, let us assume that all the $\alpha^i$'s and the $\beta^i$'s are the same. Then 
$\rho^i < \rho^j$ iff $ d^i > d^j $ and that $\rho$ will be close to a constant if the degrees are almost the same.
In other words, if the graph is ``almost" homogeneous, $L(\rho)$ would be small and hence $\rho$ will be close to the true steady state 
$\bar{p}$. We can set up an iterative procedure starting with the initial guess $p(0) = \rho$ and iterating using the Laplacian:
We define recursively, $ (p(k+1))^i$ to be the solution $\in (0 , 1)$ of the equation: 
$$R^i(p(k+1)^i) = \hat{\gamma} \, (p(k))^i \, L(p(k))^i +  \gamma_{10} (L(p(k))^i$$
This will converge rapidly to the steady state if the graph is ``almost" homogeneous.

\bigskip

\subsection{Some special cases}

\subsubsection{The homogeneous case}
If all nodes have the same matrix $Q$, the same degree $d$ and all the  non-zero weights are equal to $1$ (i.e. the network is a regular graph),
 then the stationary probability $( \bar{p} , \bar{q} ) $ is the same for all nodes and since the Laplacian vanishes on constant functions we get: 
\bear
\hat{\gamma}\,d\, \bar{p}^2  - (\alpha + \beta + \hat{\gamma}\,d) \, \bar{p} + \beta &=& 0 
\eear
\medskip
\noi
This quadratic equation has exactly one zero in the interior of $[0 , 1]$,  provided $\alpha > 0$ and $\beta > 0$. If $\hat{\gamma} = 0 $
then $\bar{p} = \frac{\beta}{\alpha + \beta} $. It is also easy to check that $\bar{p} < \frac{\beta}{\alpha + \beta} \,$ if $\hat{\gamma} >0$ 
and $\bar{p} > \frac{\beta}{\alpha + \beta} \,$if $\hat{\gamma} < 0$, so the probability strictly changes if the Markov chains are linked by a network. 
In fact, $\bar{p} \rightarrow 0$ as $\hat{\gamma} \rightarrow + \infty$  and $\bar{p} \rightarrow 1$ as $\hat{\gamma} \rightarrow - \infty$. 
Note also that even if $\alpha = 0$ there is a solution $ \bar{p} = \frac{\beta}{\hat{\gamma}} \in (0,1)$ provided $ 0 < \beta < \hat{\gamma}$ 
and if $\beta = 0$ there is a solution $\bar{p} =  1 + \frac{\alpha}{\hat{\gamma}} \in (0,1)$ provided $ 0 < \alpha < - \, \hat{\gamma}$.
On the hypercube $\{ 0 ,1 \}^N$, the probabilities are then binomially distributed. The probability 
at a state with $ k \, |0 \rangle $'s and $ l \, |1 \rangle$'s is $\bar{p}^k \, \bar{q}^l$. \\
The proof of the uniqueness and stability of the steady state in the homogeneous case can be simplified using another useful little fact about the Laplacian which we would like to record here
(the proof is elementary). 

\begin{Lemma}
 If $x^i \in (0,1)$ for all $i = 1,\ldots, N$, then
$$ \sum_{i,j}  \frac {x^i} {1 - x^i} L^i_j x^j  \geq 0  $$ 
and equality holds iff $x^i = x^j$ for all $i, j$.
\end{Lemma}

\medskip

\subsubsection {SIS model }
This is a simple epidemiological model (see \cite{ne}), corresponding to $\alpha = 0$ and $\gamma_{10} = 0 $ in our notation. {\it In the epidemiological literature, what we call $\be$ is $\gamma$, what we call $\gamma_{01}$ is $\be$, the state $|0 \rg$ is called $S$(susceptible), $|1 \rg$ is $I$(infected)}. Let us also assume that all the Markov chains are identical, so all the $\beta$'s are the same. $S$ is an absorbing steady state at each site in the absence of connections. If the network is homogeneous (a regular graph) where every node has the same degree $d$, there is another stable steady state solution (endemic equilibrium) $ \bar{p} = \frac{\beta}{ d \gamma} \in (0,1)$ provided $ 0 < \beta < d \gamma$. In the case of a general network, our proof shows that if there is an endemic equilibrium in the interior (this is true in many cases), it will be unique and stable.

\subsubsection{The case $\hat{\gamma} = 0$}
\bigskip
If we assume that the two interaction strengths are the same $\gamma_{0 1} = \gamma_{1 0} = \gamma $ and all the $\alpha$'s and $\beta$'s are equal (but we do not assume that the network is homogeneous), then the equation for the equilibrium state simplifies to:
\bear
  - (\alpha + \beta) \, \bar{p}^i  + \beta =   \gamma \sum_j L_j^i \bar{p}^j
 \eear
Averaging over $i$, we see that $\bar{p}$ is a constant equal to $\frac{\beta}{\alpha + \beta}$, so the network has no effect in this case and the ``synchronization" is perfect.

\end{document}